\documentclass[12pt]{amsart}
\usepackage{latexsym}
\usepackage{amsthm}
\usepackage{amssymb}
\usepackage{amsfonts, amsmath}
\usepackage{enumerate}

\newtheorem{Prop}{Proposition}[section]
\newtheorem{Def}[Prop]{Definition}

\newtheorem{Thm}[Prop]{Theorem}
\newtheorem{Cor}[Prop]{Corollary}
\newtheorem*{Rem}{Remark}

\title{S-rings of Krull type}
\author{\bf M. Mou\c{c}ouf}
\date{}

\begin{document}
 \maketitle
\begin{center}
{\footnotesize Department of Mathematics, Faculty of Science, Chouaib Doukkali University, Morocco\\
Email: moucouf@hotmail.com}
\end{center}
\begin{abstract}
We define and give some properties and characterizations
of S-rings of Krull type. We also determine, in the case of
an independent S-ring of Krull type $A$, the injective
dimension of the quotient category $Mod(A)/\mathcal{M}_{0}$, where
$\mathcal{M}_{0}$ is the thick subcategory of the modules with
trivial maps into the codivisorial modules.
\end{abstract}
\maketitle
 \section*{Introduction}
 Let $A$ be a ring of Krull type and $\Omega$,
(resp. $\Omega^{t}$) a defining (resp. the thin defining) family
for $A$ (see ~\cite{Mou} and ~\cite{Mo}), let $Mod(A)$ denote the
category of all unitary $A$-modules and let $\mathcal{M}_{0}$
denote the full subcategory of $Mod(A)$ consisting of all modules
$M$ such that $M_{\omega} = 0$ for all $\omega\in \Omega$, and let
$\mathcal{M}$ denote the full subcategory of all $A$-modules $M$
such that $M$ has no subobject other then $0$ belonging to
$\mathcal{M}_0$ (the elements of $\mathcal{M}$ are said to be
codivisorial).
\\\\\hspace*{1pc}In~\cite{Mou} we have proved that the
category $\mathcal{M}$ is closed under injective envelopes. So, it
is natural to wonder for which rings of Krull type the category
$\mathcal{M}_{0}$ is also closed under injective envelopes.
\\ We define S-ring of Krull type as a ring of Krull type such
that the category $\mathcal{M}_{0}$ is closed under injective
envelopes.
\\\\\hspace*{1pc}The purpose of this paper is to study
some properties of a S-ring of Krull type.
\\\hspace*{1pc}In the first section of this paper we give some properties and characterizations of
S-rings of Krull type by using c-ideals and injective
envelopes notions. We recall that in~\cite{Mo}, we have proved that
the property `` $E(M)_{\omega}\simeq E(M_{\omega})$ for all
$\omega\in \Omega^{t}$ '' is true for all codivisorial $A$-module
$M$ exactly if $A$ is an independent ring of Krull type (see
~\cite{Mo} Prop. 1.9.); in this section, we prove that this
property is true for all $A$-module $M$ exactly if $A$ is an
independent S-ring of Krull type.
\\\hspace*{1pc}In the second section we study injective dimension
in the category quotient $Mod(A)/\mathcal{M}_{0}$. Essentially, we
prove that in the case of an independent S-ring of Krull
type $A$ we have,
$$inj.dim_{A}(M) = inj.dim_{Mod(A)/\mathcal{M}_{0}}(TM)\;\text{for all}
\;\Re_{\omega}\text{-module}\; M \;\text{and}\;\omega\in \Omega^{t}.$$ and
$$inj.dim(Mod(A)/\mathcal{M}_{0}) = \sup_{\omega\in
\Omega^{t}}gl.dim(\Re_{\omega}).$$
\hspace*{1pc}We recall that in an abelian category $\mathcal{A}$ we have the followings
definitions and results:
\begin{enumerate}[-]
\item A non-empty full
subcategory $\mathcal{C}$ of $\mathcal{A}$ is called thick if the
following holds: If
$$0\longrightarrow M'\longrightarrow M\longrightarrow M''\longrightarrow 0$$ is a short exact sequence in $\mathcal{A}$, then $M$ is an object of
$\mathcal{C}$ if and only if $M'$ and $M''$ are objects of
$\mathcal{C}$.
\item It is not difficult to show that
$\mathcal{M}_0$ is a thick subcategory of  $Mod(A)$, and we can
then consider the quotient category $Mod(A)/\mathcal{M}_{0}$ of
$Mod(A)$ by $\mathcal{M}_{0}$, and the canonical functor $T:
Mod(A)\longrightarrow Mod(A)/\mathcal{M}_{0}$ (for more details, see~(\cite{Gab}, Chap. III)).
\item Throughout this paper we will use the notation and terminology of~\cite{Gri.1}, ~\cite{Gri.2}, ~\cite{Mou}, and ~\cite{Mo} for valuation ring, defining family, thin defining family and the center of a valuation.
\end{enumerate}
For the sake of completeness we give here the central definitions and notations.
\\Let $A$ be an integral domain with identity having quotient field $K$, $\Omega$ be a set of valuations of $K$, and consider the following conditions on $\Omega$
\begin{enumerate}[(1)]
\item $A=\displaystyle\bigcap_{\omega\in \Omega}\Re_\omega$,
where $\Re_\omega$ is the ring of $\omega$.
\item $\omega$ is essential for $A$; ie, $\Re_\omega$ is a quotient ring for $A$.
\item For every non-zero element $x\in K$, the set
$\{\omega\in \Omega :\; \omega(x)\neq 0\}$ is finite.
\item The valuations of $\Omega$ are pairwise independent.
\item Each $\omega$ of $\Omega$ is rank one discrete.
\end{enumerate}
Following Griffin~\cite{Gri.1} we say that $A$ is a ring of Krull type (respectively an independent ring of Krull type) if $\Omega$ satisfies $(1)$, $(2)$ and $(3)$ (respectively $(1)-(4)$). A domain with a familly $\Omega$ satisfying $(1)$-$(5)$ is said to be a krull domain.
\\In any of the above cases, the set $\Omega$ is called a defining family for $A$. Without loss of generality, we may assume that $\Omega$ does not contain two equivalent valuations.
\\A thin defining family for $A$ is a defining family $\Omega^{t}$ for $A$ such that for all $\omega\neq \upsilon \in \Omega^{t}$, we have $\Re_\omega\not\subset \Re_\upsilon$; which is equivalent to
$\wp(\omega)\not\subset \wp(\upsilon)$ for all $\omega\neq
\upsilon\in \Omega^{t}$ (see ~\cite{Mou}). The existence of a thin defining family
is confirmed by Griffin in ~\cite{Gri.2} (see Lemma 18). It is
easily seen, by using Proposition $12$ of~\cite{Gri.1}, that such
family is unique.
\\For a non-zero element $a$ of $A$ and a non-zero ideal $\vartheta$ of $A$, we denote by $\Omega(a)$ and $\Omega(\vartheta)$ the following finite subsets of $\Omega$
$$\Omega(a)=\{\omega\in \Omega/a\in \wp(\omega)\}\quad\text{and}\quad\Omega(\vartheta)=\{\omega\in \Omega/\vartheta\subseteq \wp(\omega)\}$$
where $\wp(\omega)=A\cap \wp$ and $\wp$ is the maximal ideal of the valuation $\omega$.
\\If $M$ is an $A$-module, we denote by $M_{\omega}$ the localization of $M$ at the prime ideal $\wp(\omega)$.
\\Throughout this paper $A$ will denote a ring of Krull type and $K$ the quotient
field of $A$.
\section{S-rings of Krull type}
Our objective in this section is to give some
properties and characterizations of a S-ring of Krull type,
and prove that $A$ is an independent S-ring of Krull type
if and only if $E(M_{\omega})\simeq E(M)_{\omega}$ for all
$A$-module $M$ and all $\omega\in \Omega^{t}$.
\\\\\hspace*{1pc}Let $M$ be an $A$-module. In what follows, we introduce the following subset of $M$
$$M^{c} = \{x\in M:\;Ann(x)_{c} = A\}.$$
It easily seen that
$M^{c}$ is the maximal submodule of $M$ belonging to $\mathcal
{M}_0$ ($M^{c} = ker(P_{M})$ where $P_{M}:M\longrightarrow
\prod_{\omega\in \Omega}M_{\omega}$ is the canonical mapping), and
that $M/M^{c}$ is a codivisorial module. We begin with a theorem.
\begin{Thm}\label{Thm 1}
Let $A$ be a ring of Krull type. Then the following
conditions are equivalent:
\begin{enumerate}[1)]
\item If $M$ is an $A$-module not belonging to $\mathcal {M}_0$, then $M$ has a
non-zero codivisorial submodule.
\item For all $A$-module $M$ and all
essential extension $N$ of $M$, if $M\in \mathcal {M}_0$, then so is
$N$.
\item $\mathcal {M}_0$ is closed under injective
envelopes.
\item $E(M^{c})= E(M)^{c}$ For all $A$-module $M$.
\item If $E$ is an injective $A$-module, then so is
$E^{c}$.
\item For all ideal $\vartheta$ of $A$, there exists an ideal $\vartheta'$ of $A$ such that $\vartheta'_{c}
= A$ and $\vartheta =
\vartheta_{c}\cap\vartheta'$.
\item For all ideal $\vartheta$ of $A$, the $A$-module $A/\vartheta_{c}$ is injected in $A/\vartheta$.
\item For all ideal $\vartheta$ of $A$ with $\vartheta_{c}\neq A$, there exists $a\in A - \vartheta_{c}$
such that $\vartheta:_{A}a = \vartheta_{c}:_{A}a$.
\item For all ideal $\vartheta$ of $A$ with $\vartheta_{c}\neq A$, there exists $a\in A - \vartheta_{c}$ such that $\vartheta:_{A}a$ is a
c-ideal.
\end{enumerate}
\end{Thm}
\proof~\\
$1)\Rightarrow 2)$ If $N\notin \mathcal {M}_0$ then there exists a
non-zero codivisorial submodule $N'$ of $N$. Since $N$ is an
essential extension of $M$, $N'\cap M$ is a non-zero codivisorial
submodule of $M$, and $M$ will be not in $\mathcal {M}_0$; a
contradiction.
\\$2)\Rightarrow 3)$ Trivial.
\\$3)\Rightarrow 4)$ We have
$E(M^{c})\in \mathcal {M}_0$, and hence $E(M^{c})\subset E(M)^{c}$
by the maximality of $E(M)^{c}$. $E(M)^{c}$ is an essential
extension of $M^{c}$ because $E(M)^{c}\cap M = M^{c}$. Thus
$E(M^{c}) = E(M)^{c}$.
\\$4)\Rightarrow 5)$ Let $E$ be an injective module.
Then $E(E^{c}) = E(E)^{c} = E^{ c}$. Hence $E^{c}$ is injective.
\\$5)\Rightarrow 6)$ Consider $E = E(A/\vartheta)$.
$E^{c}$ is an injective submodule of $E$, then there exists a
submodule $E'$ of $E$ such that $E'\oplus E^{c} = E$. It is easily
seen that $E'$ is codivisorial. Let $\overline{1} = x_{1} + x_{2}$
where $x_{1}\in E'$ and $x_{2}\in E^{c}$. Then $\vartheta =
Ann(x_{1})\cap Ann(x_{2})$. Since $E'$ is codivisorial,
$Ann(x_{1})$ is a c-ideal of $A$ (see ~\cite{Mou}, Proposition 1.5) contains $\vartheta$, then it also contains $\vartheta_{c}$,
and we can then replace $Ann(x_{1})$ by $\vartheta_{c}$. On the
other hand, $E^{c}\in \mathcal {M}_0$, then $Ann(x_{2})_{c} =
A$.
\\$6)\Rightarrow 7)$ We have $\vartheta =
\vartheta_c\cap \vartheta'$ where $\vartheta'_c = A$. Obviously
$\displaystyle\vartheta'\not\subset \bigcup_{\omega\in
\Omega(\vartheta)}\wp(\omega)$, for if an ideal contained in a
finite union of prime ideals then it contained in one of them. Let
$a\in \vartheta' - \displaystyle\bigcup_{\omega\in \Omega(\vartheta)}\wp(\omega)$.
Then $\vartheta:_{A}a = (\vartheta_{c}:a)\cap (\vartheta':_{A}a) =
\vartheta_{c}:a$. But since $a$ is a unit in $\Re_{\omega}$ for
all $\omega\in \Omega(\vartheta)$ it follows that $\vartheta_{c}:a
= \vartheta_{c}$. Hence $\vartheta:_{A}a = \vartheta_{c}$.
\\$7)\Rightarrow 8)$ If $\vartheta = \vartheta_{c}$,
 then we take $a = 1$. Suppose that $\vartheta\neq \vartheta_{c}$ and let $a\in A$ such that
$\vartheta:_{A}a = \vartheta_{c}$, then $(\vartheta:_{A}a)_{c} =
\vartheta_{c}:_{A}a = \vartheta_{c}$, it follows that $a\notin
\vartheta_{c}$ and that $\vartheta:_{A}a = \vartheta_{c}:_{A}a$.
\\$8)\Rightarrow 9)$ Trivial.
\\$9)\Rightarrow 1)$ Suppose that $M^{c}\subsetneq M$ and let $x\in M - M^{c}$. Then
$Ann(x)_{c}\neq A$. Hence there exists $a\in A - Ann(x)_c$ such
that $Ann(x):_{A}a$ is a c-ideal. This implies $Ann(ax)$ is a
c-ideal. This show that $Aax$ is a non-zero codivisorial submodule
of $M$ (see ~\cite{Mou}, Proposition 1.5).\endproof
\begin{Def}
A domain $A$ is said to be a S-ring of Krull type if $A$ is a ring of
Krull type and satisfies the equivalent conditions of Theorem~\ref{Thm 1}. $A$ is said to be an independent S-ring of Krull type
if $A$ is a S-ring of Krull type and an independent ring of
Krull type.
\end{Def}
The trivial examples of S-rings of Krull type are, finite
intersection of rings of valuations of a field $K$, h-local
Pr\"ufer rings, and generally Pr\"ufer rings of Krull type; since
in these cases we have $\mathcal {M}_0 = 0$.
\begin{Def}
We say that $\vartheta'$ is c-associated to $\vartheta$
if $\vartheta = \vartheta_{c}\cap \vartheta'$ and $\vartheta'_{c}
= A$.
\end{Def}
\begin{Prop}
Let $A$ be a ring of Krull type. Then $A$ is a S-ring of Krull
type if and only if $E(A/\vartheta_{c})$ is injected in
$E(A/\vartheta)$ for all ideal $\vartheta$ of $A$.
\end{Prop}
\proof
We use the conditions (7) and (9) of Theorem~\ref{Thm 1}.
\endproof
\begin{Prop}
Let $A$ be a ring of Krull type. Then we have:
\begin{enumerate}[1)]
\item Suppose that $\vartheta:_{A}a =\vartheta_{c}$, then $\vartheta_{c}\cap (aA +
\vartheta) = \vartheta$.
\item Let $\vartheta'$ be an
ideal c-associated to $\vartheta$, then there exists $p\in
\vartheta'$ such that $pA + \vartheta$ is an ideal c-associated to
$\vartheta$.
\end{enumerate}
\end{Prop}
\proof~
\\1) Notice that the assumption that $\vartheta:_{A}a =\vartheta_{c}$ ensures that $\vartheta:_{A}a$ is a c-ideal of $A$ and then $\vartheta:_{A}a=(\vartheta:_{A}a)_{c}=\vartheta_{c}:_{A}a$ (see~\cite{Mou}, Proposition 1.2 (d)). Let $d = ac + b,\;c\in A,\;b\in
\vartheta$ and suppose that $d\in \vartheta_{c}$, then $ac\in
\vartheta_{c}$, hence $c\in \vartheta_{c}:_{A}a=\vartheta_{c}$, so $d\in
\vartheta$.
\\\hspace*{1pc}2) Let $p\in \vartheta' -
\displaystyle\bigcup_{\omega\in \Omega(\vartheta)}\wp(\omega)$, it is clear that
$\vartheta = \vartheta_{c}\cap (pA + \vartheta)$, and since
$\Omega(p)\cap \Omega(\vartheta) = \emptyset$, we have $(pA +
\vartheta)_{c} = A$.
\endproof
\begin{Prop}\label{Prop 6}
Let $A$ be a ring of Krull type. The following
conditions are equivalent:
\begin{enumerate}[1)]
\item $A$ is a S-ring
of Krull type.
\item $(\vartheta:\vartheta_{c})_{c} =
A$ for all ideal $\vartheta$ of $A$.
\item $(\vartheta:\vartheta')_{c} = \vartheta_{c}:\vartheta'_{c}$ for
all ideal $\vartheta,\vartheta'$ of $A$.
\item $(\vartheta:\vartheta')_{c} = \vartheta_{c}:\vartheta'$ for all
ideal $\vartheta$ and c-ideal $\vartheta'$ of $A$.
\end{enumerate}
\end{Prop}
\proof~
\\$1)\Rightarrow 2)$. Let $t\in A$ such that $tA +
\vartheta$ is an ideal c-associated to $\vartheta$. It is clear
that $t\notin \displaystyle\bigcup_{\omega\in \Omega(\vartheta)}\wp(\omega)$. Let
$b\in \vartheta_{c}$, then $tb\in \vartheta_{c}\cap (\vartheta +
tA) = \vartheta$, hence $t\in \vartheta:\vartheta_{c}$, but
$\Omega(\vartheta:\vartheta_{c})\subset \Omega(\vartheta)$. So
$\Omega(\vartheta:\vartheta_{c}) = \emptyset$; i.e.
$(\vartheta:\vartheta_{c})_{c} =
A$.
\\$2)\Rightarrow 3)$. Let $a\in \vartheta:\vartheta'$, then $a\vartheta'\subset \vartheta$, hence
$a(\vartheta')_{c} = (a\vartheta')_{c}\subset \vartheta_{c}$; i.e.
$a\in \vartheta_{c}:\vartheta'_{c}$. Therefore
$\vartheta:\vartheta'\subset \vartheta_{c}:\vartheta'_{c}$, but
since $\vartheta_{c}:\vartheta'_{c}$ is a c-ideal (see~\cite{Mou}, Proposition 1.2), we have
$(\vartheta:\vartheta')_{c}\subset \vartheta_{c}:\vartheta'_{c}$.
Conversely, since $(\vartheta:\vartheta_{c})_{c} = A$, we have for
all $\omega\in \Omega$ there exists $t_{\omega}\in
(\vartheta:\vartheta_{c})\cap (A - \wp(\omega))$. Let $a\in
\vartheta_{c}:\vartheta'_{c}$, then $a\vartheta'_{c}\subset
\vartheta_{c}$, hence $t_{\omega}a\vartheta'_{c}\subset
t_{\omega}\vartheta_{c}\subset \vartheta$, and then
$t_{\omega}a\vartheta'\subset \vartheta$; i.e. $t_{\omega}a\in
\vartheta: \vartheta'$. Thus $a\in (\vartheta: \vartheta')_{c}$,
and $(\vartheta: \vartheta')_{c} =
\vartheta_{c}:\vartheta'_{c}$.
\\$3)\Rightarrow 4)$ is trivial.
\\$4)\Rightarrow 1)$. Let
$\vartheta$ be an ideal of $A$, then $(\vartheta:\vartheta_{c})_{c} = \vartheta_{c}: \vartheta_{c} = A$. Then there
exists $t\in (\vartheta: \vartheta_{c})\cap (A - \displaystyle\bigcup_{\omega\in \Omega(\vartheta)}\wp(\omega))$. Consider $\vartheta' = \vartheta + tA$. it is clear that $\vartheta'_{c} = A$. Let $b\in \vartheta_{c}\cap \vartheta'$, then $b = a + tc\in \vartheta_{c}$
where $a\in \vartheta,\;c\in A$, hence $tc\in \vartheta_{c}$,
therefore $c\in \vartheta_{c}$, but $t\in \vartheta:\vartheta_{c}$, then $tc\in \vartheta$. So $b\in \vartheta$. Thus
$\vartheta = \vartheta_{c}\cap \vartheta'$, where $\vartheta'_{c}= A$. By Theorem~\ref{Thm 1}, we have $A$ is a S-ring of Krull
type.
\endproof
\begin{Rem}
We note that in a any
ring of Krull type we have the following
results:
\\\hspace*{1pc}(2) is true if $\vartheta_{c} = \vartheta$,
$\vartheta_{c} = A$ or $\vartheta_{c}$ is finitely generated (see~\cite{Mou}, Proposition 1.2. (d)).
\\\hspace*{1pc}(3) is true if $\vartheta_{c} = \vartheta$, $\vartheta_{c} = A$, $\vartheta'_{c}= A$ or $\vartheta'$ is finitely generated.
\\\hspace*{1pc}(4) is true if $\vartheta_{c} = \vartheta$, $\vartheta_{c} = A$ or
$\vartheta'$ is finitely generated.
\end{Rem}
\begin{Prop} Let $A$ be a ring of Krull type. Then $A$ is a
S-ring of Krull type if and only if $E(M/M^{c})\simeq
E(M)/E(M)^{c}$ for all $A$-module $M$.
\end{Prop}
\proof
Suppose that $A$ is a
S-ring of Krull type. We have the homomorphism
$$M/M^{c}\longrightarrow E(M)/E(M)^{c}$$
$$\widetilde{m}\longrightarrow \overline{m},$$
but since $E(M)^{c}\cap M = M^{c}$, this homomorphism is
injective. On the other hand, since $E(M)^{c} = E(M^{c})$ is
injective, there is a submodule $E$ of $E(M)$ such that
$E(M)^{c}\oplus E = E(M)$. Let $x\in E(M), x\notin E(M)^{c}$, then
there exists $x_{1}\in E(M)^{c}, 0\neq x_{2}\in E$ such that $x =x_{1} + x_{2}$, and there exists $a\in A,\;0\neq m\in M$ such that
$ax_{2} = m$. It is clear that $\overline{m}\neq 0$ and that
$a\overline{x} = \overline{m}$, so the monomorphism
$$M/M^{c}\longrightarrow E(M)/E(M)^{c}$$
$$\widetilde{m}\longrightarrow \overline{m},$$ is an essential
extension. But $E(M)/E(M)^{c}\simeq E$; i.e. $E(M)/E(M)^{c}$ is
injective. Thus $E(M)/E(M)^{c}\simeq E(M/M^{c})$. Conversely,
suppose that $E(M/M^{c})\simeq E(M)/E(M)^{c}$ for all $A$-module $M$. let $M\in \mathcal {M}_0$, then $M = M^{c}$, hence
$E(M)/E(M)^{c}\\ = 0$. Therefore $E(M) = E(M)^{c}$ for all $A$-module $M\in \mathcal {M}_0$. $\mathcal {M}_0$ is then closed under injective envelopes. Thus $A$ is a S-ring of Krull
type.
\endproof
The next proposition gives a characterization of an independent
S-ring of Krull type.
\begin{Prop}\label{Prop 8}
Let $A$ be a ring of Krull type. Then $A$ is an independent S-ring of Krull
 type if and only if $E(M)_{\omega}\simeq E(M_{\omega})$ for all
$A$-modules $M$ and all $\omega\in \Omega^{t}$. In this case we
have, $E_{\omega}$ is an injective $A$-module for all injective
$A$-module $E$ and all $\omega\in \Omega^{t}$.
\end{Prop}
\proof
The ``if part'' follows immediately from 3) of Theorem~\ref{Thm 1}
and Proposition 1.9. of ~\cite{Mo}. For the ``only if part'', we have $E(M)/E(M)^{c}\simeq E(M/M^{c})$, then
$E(M)_{\omega}/E(M)^{c}_{\omega}\simeq E(M/M^{c})_{\omega}$, hence $E(M)_{\omega}\simeq E(M/M^{c})_{\omega}$, but $M/M^{c}$ is codivisorial and $A$ is an independent ring of Krull type, then by Proposition 1.9. of ~\cite{Mo} we have $E(M/M^{c})_{\omega}\simeq
E((M/M^{c})_{\omega})=E(M_{\omega}/M^{c}_{\omega})=E(M_{\omega})$, so $E(M)_{\omega}\simeq E(M_{\omega})$. The last assertion is trivial.
\endproof
Let $A$ be a ring of Krull type and $\Omega$ be a defining family for $A$. Let us denote by $\xi$ the set of all ideals $\vartheta$ of $A$ such that
$\vartheta\Re_{\omega}$ is principal for all $\omega\in \Omega$. Then we have the following result.
\begin{Prop}\label{Prop 1.10}
Let $A$ be an independent ring of Krull type and let $\vartheta$ be an ideal of $A$ that lies in $\xi$. Then
there exists $t\in A$ such that $\vartheta_{c} =\vartheta:t$.
\end{Prop}
\proof If $\vartheta_{c}=A$, then $t=0$ does the job, so let us assume that $\vartheta_{c}\neq A$; i.e, $\Omega^{t}(\vartheta)\neq \emptyset$.
We have $\vartheta_{c} =\displaystyle
\bigcap_{\omega\in \Omega^{t}(\vartheta)}(a_{\omega}\Re_{\omega}\cap
A)$ where $a_{\omega}\in \vartheta$ for all $\omega\in
\Omega^{t}(\vartheta)$. By Lemma 2.13. of~\cite{Mou}, there
exists $a\in \vartheta$ such that $\omega(a) = \omega(a_{\omega})$
for all $\omega\in \Omega^{t}(\vartheta)$, and By the
approximation theorem for independent ring of Krull type (see~\cite{Gri.2} and \cite{Knebusch}), there exists $t\in A$ such that $\omega(t) = 0$ for all
$\omega\in \Omega^{t}(\vartheta)$ and $\omega(t) = \omega(a)$ for
all $\omega\in \Omega^{t}(a) - \Omega^{t}(\vartheta)$. Since
$t\notin \displaystyle\bigcup_{\Omega^{t}(\vartheta)}\wp(\omega)$, we have
$\vartheta:t\subset \vartheta_{c}$. Conversely, let $d\in
\vartheta_{c}$, then
$$\omega(td) \left\{\begin{array}{ll}
=\omega(d)\geq \omega(a) &\text{ if \;\;} \omega\in\Omega^{t}(\vartheta) \\
\geq \omega(t) = \omega(a) & \text{ if \;\;}\omega\in\Omega^{t}(a)
- \Omega^{t}(\vartheta) \\\geq 0 = \omega(a) & \text{ if
\;\;}\omega\notin\Omega^{t}(a)
\end{array}\right.$$ Hence $td \in aA\subset \vartheta$. Thus
$\vartheta:t = \vartheta_{c}$
\endproof
By using Proposition~\ref{Prop 1.10} and the fact that in the case of Krull domains, $\xi$ is the set of all ideals, we find the following result due to Beck:
\begin{Cor} (Beck, 1971, Corollary 3.5)
 A Krull domain is a S-ring of Krull type.
\end{Cor}
Recall~\cite{Mo} that a ring of Krull type $A$ is said to be a $(P)$-ring if $A$ possesses a defining family consisting of valuations with principal maximal ideals.
\begin{Prop}
Let $A$ be a $(P)$-ring. Then we have
\begin{enumerate}[1)]
\item $Ass(A/\vartheta_{c})=\Omega^{t}(\vartheta)$ for every ideal $\vartheta\in \xi$.
\item If $A$ is an independent ring of Krull type, then $\Omega^{t}(\vartheta) = \{\wp\in Ass(A/\vartheta):\;\wp = \wp_{c}\}$ for every ideal $\vartheta\in \xi$.
\item If $A$ is a S-ring of Krull type, then $\Omega^{t}(\vartheta) = \{\wp\in Ass(A/\vartheta):\;\wp = \wp_{c}\}$ for every ideal of $A$.
\end{enumerate}
\end{Prop}
\proof~
\begin{enumerate}[1)]
\item In view of Corollary 1.6 (3) and Lemma 1.1 of~\cite{Mo}, we have $E(A/\vartheta_{c})\simeq
\oplus_{\omega\in
\Omega^{t}(\vartheta)}E(A/\wp(\omega))$, therefore $Ass(A/\vartheta_{c})=\Omega^{t}(\vartheta)$.
 \item Since $A$ is an independent ring of Krull type, it follows from~\ref{Prop 1.10} that $Ass(A/\vartheta_{c})\subseteq Ass(A/\vartheta)$ and then $\Omega^{t}(\vartheta)\subseteq Ass(A/\vartheta)$. Now suppose that $\wp$ is prime c-ideal of $A$ that lies in $Ass(A/\vartheta)$. Then there exists $s\in A$ such that $\wp=\vartheta : s$. Hence $\wp=(\vartheta : s)_{c}=\vartheta_{c} : s$, so $\wp\in \Omega^{t}(\vartheta)$.
\item Follows from Theorem~\ref{Thm 1} (7) by using the same proof as in 2).
\end{enumerate}
\endproof
We conclude this section with three corollaries of the above
results.
\begin{Cor}(Beck, 1971, Corollary 3.2)
Let $A$ be a Krull domain and Let $\vartheta\neq 0$ be an
ideal of $A$. Let $\wp$ be a prime minimal ideal of $A$ containing $\vartheta$. Then $\wp\in Ass(A/\vartheta)$.
\end{Cor}
\begin{Cor}
Let $A$ be a Krull domain and let
$\vartheta,\vartheta'$ be fractional ideals of $A$. Then
$(\vartheta:_{k}\vartheta')_{c} =
\vartheta_{c}:_{k}\vartheta'_{c}$.
\end{Cor}
\proof
Let $r,s\in A$
such that $r\vartheta\subset A$ and $s\vartheta'\subset A$. Since
$A$ is completely integrally closed, we have
$r\vartheta:_{k}s\vartheta' = r\vartheta: s\vartheta'$, but by
Proposition~\ref{Prop 6}, we have $(r\vartheta: s\vartheta')_{c} =(r\vartheta)_{c}: (s\vartheta')_{c} = r\vartheta_{c}:s\vartheta'_{c}$, then $(r\vartheta:_{k} s\vartheta')_{c} =r\vartheta_{c}:_{k} s\vartheta'_{c}$, but $(r\vartheta:_{k}s\vartheta')_{c} =(\vartheta:_{k}\vartheta':_{k}sr^{-1})_{c}=(\vartheta:_{k}\vartheta')_{c}:_{k}sr^{-1})$ and
$r\vartheta_{c}:_{k} s\vartheta'_{c} = (\vartheta_{c}:_{k}\vartheta'_{c}):_{k}sr^{-1}$. Thus $(\vartheta:_{k}\vartheta')_{c}=\vartheta_{c}:_{k}\vartheta'_{c}$.
\endproof
\begin{Cor} Let $A$ be a S-ring of Krull type. Let $\vartheta$ be an irreducible ideal of $A$ such that $\vartheta_{c}\neq A$. Then $\vartheta$ is a c-ideal.
\end{Cor}
\section{The Injective Dimension Of The Quotient Category $Mod(A)/\mathcal{M}_{0}$}
In this section, we prove that over an independent S-ring of Krull type, the injective dimension
of the quotient category $Mod(A)/\mathcal{M}_{0}$ is exactly equal to $\sup_{\omega\in \Omega^{t}}inj.dim(\Re_{\omega})$.
\\ We note that:
\begin{enumerate}[-]
\item The thick subcategory $\mathcal{M}_{0}$ of $Mod(A)$ is localizing, since for
all $A$-module $M$: among the subobjects of $M$ belonging to $\mathcal{M}_{0}$ there is a maximal one (see Corollary 1, p. 375 of ~\cite{Gab}).
\item $Mod(A)/\mathcal{M}_{0}$ is a category with injective
envelopes, and every injective object of this category is
isomorphic to an object $TI$ where $I$ is a codivisorial injective
$A$-module (see Corollary 2, p. 375 of~\cite{Gab}).
\item If $A$ is a S-ring of Krull type, then $TJ$ is an
injective object of $Mod(A)/\mathcal{M}_{0}$ for every injective
$A$-module $J$ (see Corollary 3, p. 375 of~\cite{Gab}).
\end{enumerate}
We begin with a proposition.
\begin{Prop}\label{Prop 21}
Let $A$ be a ring of Krull type, $\Omega$ a defining family for $A$, $M,N$ two any
$A$-modules, and let $f$ be any $A$-homomorphism from $M$ to $N$.
Then:
\begin{enumerate}[1)]
\item Suppose that
$f_{\omega}:M_{\omega}\longrightarrow N_{\omega}$ is an essential
extension of $M_{\omega}$ for all $\omega\in \Omega$. Then the
homomorphism
$$\overline{f}:M/M^{c}\longrightarrow N/N^{c}$$
$$\overline{m}\longrightarrow \overline{f(m)}$$ is an
essential extension of $M/M^{c}$.
\item Conversely, suppose that $f$ is an essential extension of $M$. Then
$f_{\omega}:M_{\omega}\longrightarrow N_{\omega}$ is an essential
extension of $M_{\omega}$ for all $\omega\in \Omega^{t}$, in the
following cases:
\\\hspace*{1pc}(i) $M$ is codivisorial and $A$ is
an independent ring of Krull type.
\\\hspace*{1pc}(ii) $A$ is an
independent S-ring of Krull type.
\end{enumerate}
\end{Prop}
\proof~
\begin{enumerate}[1)]
\item It is
easily seen that $f(M^{c})\subset N^{c}$, and then $\overline{f}$
is well defined, but since $f_{\omega}$ is injective for all
$\omega\in \Omega$, it follows easily that $\overline{f}$ is also
injective. On the other hand, let $\overline{H}$ be a submodule of
$N/N^{c}$ ($N^{c}\subset H$) such that $\overline{H}\cap
\overline{f}(M/M^{c}) = 0$. Then $\overline{H}_{\omega}\cap
(\overline{f}(M/M^{c}))_{\omega} = 0$, but since $N^{c}_{\omega} =
0$, it follows that $H_{\omega}\cap f(M)_{\omega} = 0$, and then
$H_{\omega} = 0$ for all $\omega\in \Omega$; i.e. $H\in
\mathcal{M}_{0}$. Therefore $H\subset N^{c}$. Thus $\overline{H} =
0$. This implies that $\overline{f}$ is an essential extension of
$M/M^{c}$.
\item The two statements follow easily from~\cite{Mo} Proposition 1.9 and from Proposition~\ref{Prop 8}, respectively.
\end{enumerate}
\endproof
\begin{Cor}
Let $A$ be an  independent S-ring of Krull type. If $M$ is
 an $\Re_{\omega}$-module for some $\omega\in
\Omega^{t}$, then $inj.dim_{\Re_{\omega}}(M) = inj.dim_{A}(M)$.
\end{Cor}
\proof
Follows immediately from Proposition~\ref{Prop 8}.
\endproof
The next proposition is a generalization of Proposition 2.10. of ~\cite{Beck}.
\begin{Prop}
Let $A$ be a ring of Krull type, and let $M$ be any
$A$-module. Then
$$inj.dim_{A}(\oplus_{\omega\in \Omega}M_{\omega})\leq \sup_{\omega\in \Omega}gl.dim(\Re_{\omega}).$$
\end{Prop}
\proof
The assertion is trivial if $\sup_{\omega\in \Omega}gl.dim(\Re_{\omega}) =\infty$. Suppose that $\sup_{\omega\in \Omega}gl.dim(\Re_{\omega}) = n$, then by considering $M_{\omega}$ as an $\Re_{\omega}$-module, we can find an exact sequence
$$0\longrightarrow M_{\omega}\longrightarrow (E_{0})_{\omega}\longrightarrow\cdots\longrightarrow (E_{n})_{\omega}\longrightarrow 0$$ where $(E_{i})_{\omega}$ is an injective $\Re_{\omega}$-module for all $i\in \{1,...,n\}$. But since $\Re_{\omega}$ is a quotient ring of $A$, we have
$(E_{i})_{\omega}$ is an injective $A$-module for all $i\in \{1,...,n\}$, and then $\oplus_{\omega\in \Omega}(E_{i})_{\omega}$
is an injective $A$-module (see ~\cite{Mou}, Corollary 2.9.).
Hence, the assertion follows from the fact that the sequence
$$0\longrightarrow \oplus_{\omega\in \Omega}M_{\omega}\longrightarrow \oplus_{\omega\in \Omega}(E_{0})_{\omega}\longrightarrow\cdots\longrightarrow \oplus_{\omega\in \Omega}(E_{n})_{\omega}\longrightarrow 0$$ is exact.
\endproof
\begin{Prop}\label{Prop 24}
Let $A$ be a ring of Krull type, and
let $M$ be an $\Re_{\omega}$-module for some $\omega\in
\Omega^{t}$. Then:
\begin{enumerate}[1)]
\item $inj.dim_{A}(M)\leq inj.dim_{Mod(A)/\mathcal{M}_{0}}(TM)$ if $A$ is an independent ring of
Krull type.
\item $inj.dim_{Mod(A)/\mathcal{M}_{0}}(TM)\leq inj.dim_{A}(M)$ if $A$ is a
S-ring of Krull type .
\item If $A$ is an independent S-ring of Krull type, then$$inj.dim_{A}(M) =inj.dim_{Mod(A)/\mathcal{M}_{0}}(TM).$$
\end{enumerate}
\end{Prop}
\proof~
\begin{enumerate}[1)]
\item If $inj.dim_{Mod(A)/\mathcal{M}_{0}}(TM)=\infty$, we have done. We may then assume $inj.dim_{Mod(A)/\mathcal{M}_{0}}(TM) = m$.
Let $$0\longrightarrow TM\longrightarrow C_{0}\longrightarrow\cdots\longrightarrow C_{m}\longrightarrow 0$$ an injective resolution of $TM$, we know
that there exists $I_{1},...,I_{m}$ injective codivisorial $A$-modules such that $C_{i}\simeq TI_{i},\;1\leq i\leq m$. On the
other hand, consider the abelian category $Mod(A)_{\omega}$ of all $\Re_{\omega}$-modules, and consider
$G =\cdot\otimes_{A}\Re_{\omega} : Mod(A)\longrightarrow Mod(A)_{\omega}$ which is an exact functor. Since $GN = 0$ for all $N\in \mathcal{M}_{0}$, there is a (unique up to isomorphism) functor $H:Mod(A)/\mathcal{M}_{0}\longrightarrow Mod(A)_{\omega}$, such
that $G = HT$ (see ~\cite{Gab}, Corollary 2. p. 368), and since $G$ is exact, $H$ is also exact (see ~\cite{Gab}, Corollary 3. p.369). Now consider the exact sequence $$0\longrightarrow TM\longrightarrow TI_{0}\longrightarrow\cdots\longrightarrow TI_{m}\longrightarrow 0$$
We apply the exact functor $H$ to this sequence, and drive the
following exact sequence
$$0\longrightarrow M_{\omega} = M\longrightarrow (I_{0})_{\omega}\longrightarrow\cdots\longrightarrow (I_{m})_{\omega}\longrightarrow 0$$ where
$(I_{0})_{\omega},...,(I_{m})_{\omega}$ are injective by Proposition 1.9. of ~\cite{Mo}.
Thus $inj.dim_{A}(M)\leq inj.dim_{Mod(A)/\mathcal{M}_{0}}(TM)$.
\item Follows from the fact that if $A$ is a S-ring of
Krull type, then $\mathcal{M}_{0}$ is closed under injective
envelopes, and $TE$ is an injective object of
$Mod(A)/\mathcal{M}_{0}$ if $E$ is an injective
$A$-module.
\item Trivial.
\end{enumerate}
\endproof
\begin{Thm}
Let $A$ be an independent S-ring of Krull type. Then
$$inj.dim(Mod(A)/\mathcal{M}_{0}) = \sup_{\omega\in \Omega^{t}}gl.dim(\Re_{\omega}).$$
\end{Thm}
\proof
Suppose in the first case that
$\sup_{\omega\in \Omega^{t}}gl.dim(\Re_{\omega}) = n$, let $M$ be
any $A$-module, and let $0\longrightarrow M\longrightarrow
E_{0}\longrightarrow E_{1}\longrightarrow\cdots$ be a minimal
injective resolution of $M$. Then, by Proposition~\ref{Prop 8} and
Proposition~\ref{Prop 21} 2), (ii), $0\longrightarrow M_{\omega}\longrightarrow (E_{0})_{\omega}\longrightarrow (E_{1})_{\omega}\longrightarrow\cdots$ is a minimal injective resolution of $M_{\omega}$ as an $\Re_{\omega}$-module, so
$(E_{i})_{\omega} = 0$ for all $i\geq n$ and for all $\omega\in \Omega^{t}$. It follows that $E_{i}\in \mathcal{M}_{0}$ for all
$i\geq n$, and then $0\longrightarrow TM\longrightarrow TE_{0}\longrightarrow\cdots\longrightarrow TE_{n}\longrightarrow 0$
is an injective resolution of the object $TM$ of $Mod(A)/\mathcal{M}_{0}$. Thus $inj.dim(Mod(A)/\mathcal{M}_{0})\leq n$. On the other hand, since
$\sup_{\omega\in \Omega^{t}}gl.dim(\Re_{\omega}) = n$, there
exists $\omega\in \Omega^{t}$ such that $gl.dim(\Re_{\omega}) =n$, hence there exists an $\Re_{\omega}$-module $M$ such that
$inj.dim_{\Re_{\omega}}(M) = n$. By the preceding proposition, we have $n\leq inj.dim(Mod(A)/\mathcal{M}_{0})$. So $inj.dim(Mod(A)/\mathcal{M}_{0}) = n$. For the case, $\sup_{\omega\in \Omega^{t}}gl.dim(\Re_{\omega}) =\infty$, we use the fact that
for all $n\in \mathbb{N}$, there exists $\omega\in \Omega^{t}$ and an $\Re_{\omega}$-module $M$ such that $inj.dim_{\Re_{\omega}}(M)\geq n$. It is sufficient to use Proposition~\ref{Prop 24}, to have $inj.dim(Mod(A)/\mathcal{M}_{0}) =\infty$. Thus $inj.dim(Mod(A)/\mathcal{M}_{0}) = \sup_{\omega\in \Omega^{t}}gl.dim(\Re_{\omega})$.
\endproof
\begin{Cor}
Let $A$ be an independent S-ring of Krull type and $K$ its quotient field.
Then we have:
\begin{enumerate}[1)]
\item $inj.dim(Mod(A)/\mathcal{M}_{0}) =
1$ if and only if $A$ is a Krull domain (we assume that $A$ is not
a field).
\item $gl.dim(A) = \sup_{\eta\in Max(A)}gl.dim(A_{\eta})$ if $A$ is a Pr\"ufer ring and $Max(A)$ is the set
of maximal ideals of $A$.
\end{enumerate}
\end{Cor}
\proof~
\begin{enumerate}[1)]
\item Follows immediately
from the fact that $gl.dim(\Re_{\omega})\leq 1$ if and only if
$\Re_{\omega}$ is either a field or a discrete valuation
domain.
\item Follows from the fact that in the case of
Pr\"ufer ring of Krull type, we have $\mathcal{M}_{0} = 0$ and
$\Omega^{t} = \{\omega: \Re_{\omega} = A_{\rho},\;\rho\in
Max(A)\}$.
\end{enumerate}
\endproof

\end{document}